\documentclass[12pt,reqno]{article}

\usepackage[usenames]{color}
\usepackage{amssymb}
\usepackage{graphicx}
\usepackage{amscd}
\usepackage{ulem}

\usepackage[colorlinks=true,
linkcolor=webgreen,
filecolor=webbrown,
citecolor=webgreen]{hyperref}

\definecolor{webgreen}{rgb}{0,.5,0}
\definecolor{webbrown}{rgb}{.6,0,0}

\usepackage{color}
\usepackage{fullpage}
\usepackage{float}

\usepackage{graphics,amsmath,amssymb}
\usepackage{amsthm}
\usepackage{amsfonts}
\usepackage{latexsym}
\usepackage{epsf}

\def\N{\mathbb N}

\begin{document}

\begin{center}
\vskip 1cm{\LARGE\bf An Improved Upper Bound for the Sum-free Subset Constant}
\vskip 1cm
\large
Mark Lewko \\
Department of Mathematics \\
University of Texas at Austin \\
USA \\
\href{mailto:mlewko@math.utexas.edu}{\tt mlewko@math.utexas.edu} \\
\end{center}

\vskip .2 in

\begin{abstract}
We show that the optimal constant in Erd\"{o}s' sum-free subset theorem cannot be larger than $11/28 \approx .393$.
\end{abstract}

\newtheorem{theorem}{Theorem}
\newtheorem{corollary}[theorem]{Corollary}
\newtheorem{lemma}[theorem]{Lemma}
\newtheorem{proposition}[theorem]{Proposition}
\newtheorem{conjecture}[theorem]{Conjecture}
\newtheorem{defin}[theorem]{Definition}
\newenvironment{definition}{\begin{defin}\normalfont\quad}{\end{defin}}
\newtheorem{examp}[theorem]{Example}
\newenvironment{example}{\begin{examp}\normalfont\quad}{\end{examp}}
\newtheorem{rema}[theorem]{Remark}
\newenvironment{remark}{\begin{rema}\normalfont\quad}{\end{rema}}

\section{Introduction}

We say a set of natural numbers $A$ is sum-free if there is no solution to the equation $x+y=z$ with $x,y,z \in A$. The following is a well-known theorem of Erd\"{o}s \cite{EO1}.

\begin{theorem}Let $A$ be a finite set of natural numbers. There exists a sum-free subset $S \subseteq A$ such that $|S| \geq \frac{1}{3}|A|$.
\end{theorem}

The proof of this theorem is a common example of the probabilistic method and appears in many textbooks. Alon and Kleitman \cite{AK1} have observed that Erd\"{o}s' argument essentially gives the theorem with the slightly stronger conclusion $|S| \geq \frac{|A|+1}{3}$. Bourgain \cite{B01} has improved this further, showing that the conclusion can be strengthened to $|S| \geq \frac{|A| + 2}{3}$. Bourgain's estimate is sharp for small sets, and improving it for larger sets seems to be a difficult problem. There has also been interest in establishing upper bounds for the problem. It seems likely that the constant $1/3$ cannot be replaced by a larger constant, however this is an open problem. In Erd\"{o}s' 1965 paper, he showed that the constant $\frac{1}{3}$ could not be replaced by a number greater than $3/7 \approx .429$ by considering the set $\{2,3,4,5,6,8,10\}$. In 1990, Alon and Kleitman \cite{AK1} improved this to $12/29 \approx .414$. In a recent survey of open problems in combinatorics \cite{CL1}, it is reported that Malouf has shown the constant cannot be greater than $4/10 = .4$.  While we have not seen Malouf's proof, we note that this can be established by considering the set $\{1,2,3,4,5,6,8,9,10,18\}$. The purpose of this note is to further improve on these results by showing that the optimal constant cannot be greater than $11/28 \approx .393$.

Let $A \subset \N$ be a finite set with largest sum-free subset of size $l$. We define the sum-free subset constant of $A$ to be $\delta(A) = l/|A|$. For an integer $d$, let $dA$ denote the dilation of $A$ by $d$, that is $dA=\{da : a \in A\}$. Notice that if $d_{1}A \cap d_{2} A = \emptyset$, then the sum-free subset constant of the set $d_{1}A \cup d_{2}A$ is at most $\delta$ where $\delta$ is the sum-free subset constant of the set $A$. Thus if we demonstrate a single set with sum-free subset constant $\delta$, there exist arbitrarily large sets with constant at most $\delta$.

Our main result is the following theorem.

\begin{theorem}The largest sum-free subset of the set
\[A:=\{1,2,3,4,5,6,7,8,9,10,11,12,13,14,15,16,17,18,20,22,24,25,26,27,30,34,50,54\}\]
contains $11$ elements.
\end{theorem}

Of course, this theorem can be computationally verified by checking that each of the ${{28} \choose {12}} = 30,421,755$ subsets of size $12$ contain a sum. Here we give a human-verifiable proof. The set $A$ was the set with the smallest sum-free subset constant found in the course of an extensive computer search. It is easy to see that $A$ does contain a sum-free subset of size $11$: just consider the set of odd elements of $A$, $\{1,3,5,7,9,11,13,15,17,25,27\}$.

We admit that our proof is longer and more tedious than we (or the
reader) might prefer. Much of the length, however, is due to our
attempt to make the proof easily readable.

\section{Strategy}

The proof proceeds by showing that any subset $S \subset A$ of cardinality $12$ contains a sum. We start by writing our set in a strategic tabular form.

\begin{center}
\begin{tabular}{|c|c|c|c|c|c|c|c|c|c|c|c|c|}
  \hline
  1 & 2 & 3 & 4 & 5 & 6 & 7 & 8 & 9 & 10 & 11 & 12 & 13 \\ \hline
  1 & 3 & 5 & 7 & 8 & 9 & 11 & 13 & 15 & 17 & 25 & 27 & 24 \\
  2 & 6 & 10 & 14 & 16 & 18 & 22 & 26 & 30 & 34 & 50 & 54 &  \\
  4 & 12 & 20 &  & &  &  &  &  &  &  &  &  \\
  \hline
\end{tabular}
\end{center}

The columns of this table represent doubling relations. For
instance, the entries in the first column are $1$, $2 = 1+1$,
and $4 = 2+2$.  However, we have not reflected all doubling
relations with this formatting. For example, $8=4+ 4$, but $8$
is not included in the first column. The choice of which relations
to reflect with this formatting was subjective and intended to aid
in the readability of the proof.

Notice that a sum-free subset of $A$ can contain at most one element
from each of the ten columns $4-13$. Also, a sum-free subset of $A$
can contain at most two elements from each of columns $1-3$, and if such a set
does contain two elements of one of these columns, these two
elements must be exactly the first and last element of the column.

In what follows, $S$ will always denote a hypothetical sum-free subset of $A$ of size $12$. In addition, we will  denote the union of the elements of columns $1-3$ as $B$ and $C=S\cap B$. Since at most one element from each of the columns $4-13$ can be in $S$, we have that $|C| \geq 2$. Alternatively, by the remarks above if $|S \cap B| = 6$, the set $S$ would have to contain exactly the first and last element of the columns $1-3$. This would imply $1,4,5 \in S$, which would contradict the fact that $S$ is sum-free.

The following sections will consider the possible cases $|C| = 2$, $|C| =3$, $|C|=4$, and $|C| =5$ separately. We will often record the sum relation justifying a particular statement in brackets following that statement. For example, if $\{1,2\} \subset S$ we might write: we have that $3 \notin S$ [$1+2=3$].

\section{Case 1: $|C| = 2$}

If $|S|=12$ and $|C| = 2$, it follows that $S$ contains exactly one element from each of the ten columns $4-13$. In particular, this implies that $24 \in S$. Also, $S$ can contain at most one element from the set $\{11,13\}$ since $11+13=24$. We consider three cases based on the possible intersection of the set $\{11,13\}$ with $C$.

\textbf{Case I: $\{11,13\}\cap S = \emptyset$}. If $\{11,13\}\cap S = \emptyset$, then $22,26 \in S$. Hence $50 \notin S$ [$26+24 =50$] and $25 \in S$. It follows that $C$ does not contain $1$ [$24+1=25$], $2$ [$22+2=24$], $3$ [$22+3=25$], $4$ [$22+4=26$], or $12$ [$12+12=24$]. We conclude that $C \subset \{6,5,10,20\}$. If we suppose that $6 \in C$, we then have that $8 \in S$ and $16 \notin S$ [$6+16=22$]; $7 \in C$ and $14 \notin C$ [$8+6=14$]. Thus $30 \notin C$ [$6+24=30$] and $15 \notin C$ [$7+8=15$]. It follows that $S$ must not contain an element of column $9$. This implies that $|S|<12$ and contradicts our hypothesis. We may now assume that $C \subset \{5,10,20\}$, hence $C=\{5,20\}$. This implies that $25 \notin S$, however this contradicts a claim above.

\textbf{Case II: $\{11,13\}\cap S = \{11\}$}. Since $\{11,13\}\cap S = \{11\}$, we have that $26 \in S$. As above, we also can assume that $24 \in S$. Also $30 \in S$ and $15 \notin S$ [$11+15=26$]; $27 \in S$ and $54 \notin S$ [30+24=54]; $8 \in S$ and $16 \notin S$ [$16+11=27$]. Hence we have $34 \notin S$ and $17 \in S$ [$8+26=34$]; $50 \notin S$ and $25 \in S$ [$26+24=50$]. However, we have reached a contradiction since $8,17,25 \in S$ and $8+17=25$.

\textbf{Case III: $\{11,13\}\cap S = \{13\}$ and $7\in S$}. If $\{11,13\}\cap S = \{13\}$ then $22 \in S$. Also, as above, we have $24 \in S$. Furthermore, $15 \notin S$ and $30 \in S$ [$7+15=22$]. Next $34 \in S$ and $17 \notin S$ [$7+17=24$]. Also, $27 \in S$ and $54 \notin S$ [$24+30=54$]. However, we have reached a contradiction since $7+27=34$.

\textbf{Case IV: $\{11,13\}\cap S = \{13\}$ and $7\notin S$}. Since $\{11,13\}\cap S = \{13\}$, we have that $22 \in S$. Also $14 \in S$ since $7 \notin S$, and, as above, we have that $24 \in S$. Thus, $9 \notin S$ and $18 \in S$ [$9+13=22$]; $16 \in S$ and $8 \notin S$ [$8+14=22$].

We now have that $1 \notin C$ [$13+1=14$], $2\notin C$ [14+2=16], $4 \notin C$ [$14+4=18$], $3 \notin C [13+3=16]$, $6 \notin C$ [$18+6=24$], $12 \notin C$ [$12+12=24$], $5 \notin C$ [$13+5=18$], $10 \notin C$ [$14+10=24$]. This implies that $|C| \leq 1$ and gives a contradiction.

\section{Case 2: $|C| = 5$}

If $|C|=5$, then by the pigeonhole principle $S$ contains two elements in two of the columns $1-3$. It is easy to see that this implies ${3,5,12,20} \in C$. Since $2 \notin C$ [$2+3=5$], we have either $C=\{1,3,5,12,20\}$ or $C=\{4,3,5,12,20\}$. We consider these cases separately. Since both of these sets contain $12$, we can conclude that $24 \notin S$ [$12+12=24$]. Thus we need to show that of the nine columns $4-12$, at least three do not contain any element of $S$.

\textbf{Case I: C=\{1,3,5,12,20\}}. As remarked above, we have that $24 \notin S$. We now have that $25 \notin S$ [$20+5=25$], $8 \notin S$ [$3+5=8$], $17 \notin S$ [$5+12=17$], $11 \notin S$ [$11+1=12$], $15 \notin S$ [$12+3=15$], $7 \notin S$ [$7+5=12$], $13 \notin S$ [$12+1=13$], $9 \notin S$ [$9+3=12$]. Let us record these results in tabular form (elements we have ruled out have a strikethrough).

\begin{center}
\begin{tabular}{|c|c|c|c|c|c|c|c|c|}
  \hline
   4 & 5 & 6 & 7 & 8 & 9 & 10 & 11 & 12  \\ \hline
  \sout{7} & \sout{8} &\rlap{|} 9 & \sout{11}  &\sout{13}  & \sout{15}  &\sout{17}  & \sout{25}  & 27  \\
    14 & 16 & 18 & 22 & 26 & 30 & 34 & 50 & 54   \\
   \hline
\end{tabular}
\end{center}

Also notice that $S$ can contain only one element from $\{14,26\}$ [$14+12=26$], only one element from $\{18, 30\}$ [$18+12 = 30$], and only two elements from $\{16, 34, 50\}$ [$16+34 = 50$]. This shows that three of the columns $4-12$ cannot contain any element of the set $S$. As remarked above, this suffices to complete the proof of this case.

\textbf{Case II: C=\{4,3,5,12,20\}}. As remarked above, we have that $24 \notin S$. Also $7 \notin S$ [$3+4=7$], $8\notin S$ [$3+5=8$], $16 \notin S$ [$4+12=16$], $9 \notin S$ [$9+3=12$], $15 \notin S$ [$3+12=15$], $17 \notin S$ [$5+12=17$], $25 \notin S$ [$5+20=25$]. We record this information in the following table.

\begin{center}
\begin{tabular}{|c|c|c|c|c|c|c|c|c|}
  \hline
   4 & 5 & 6 & 7 & 8 & 9 & 10 & 11 & 12  \\ \hline
  \sout{7} &\sout{8} &\sout{9} &  11 &  13 &\sout{15} &\sout{17} &\sout{25} & 27  \\
    14 & \sout{16} & 18 & 22 & 26 & 30 & 34 & 50 & 54   \\
   \hline
\end{tabular}
\end{center}

Now we note that at most 1 element of the following two disjoint sets can be contained in $S$: $\{14,18\}$ [$4+14=18$], $\{30,34\}$ [$4+30=34$]. However, combining this with the fact that no element of column $5$ is contained in $S$, we conclude there are at least three columns among $4-12$ that do not contain elements of $S$. This completes the proof.

\section{Case 3: $|C| = 4$}

If $|C|=4$, we must have that at least one of the columns $1-3$ contains two elements of $S$. Thus either $\{1,4\} \subset S$, $\{3,12\} \subset S$, or $\{5,20\} \subset S$. We consider these cases separately. Since we suppose that $|C|=4$, we must show that three of the ten columns $4-13$ do not contain elements of $S$.

\textbf{Case I: $\{1,4\} \subset S$}. Let us record the implications of the hypothesis $\{1,4\} \subset S$ in the following table.

\begin{center}
\begin{tabular}{|c|c|c|}
  \hline
  1 &  2 &  3  \\ \hline
  1 & \sout{3} & \sout{5}  \\
  \sout{2} & 6 & 10   \\
  4 & 12 & 20   \\
  \hline
\end{tabular}
\end{center}

We may conclude that $C$ is one of the following sets: $\{1,4,6,20\}$, $\{1,4,12,10\}$, $\{1,4,12,20\}$.

We start by assuming that $C = \{1,4,6,20\}$. First we notice that $24 \notin S$ [$4+20=24$], and thus it suffices to show that $2$ of the columns $4-12$ contain no elements of $S$. Next, we have that $8 \notin S$ [$4+4=8$] and $16 \notin S$ [$4+16=20$]. It follows that no element of column $5$ is in $S$, and hence we may assume that one element of each of the remaining columns (among $4-12$) contains an element of $S$. Thus $14 \in S$ and $7 \notin S$ [$1+6=7$]. However, this gives a contradiction since we now have that $6,14$ and $20$ are in $S$ [$6+14=20$].

We now suppose that $C=\{1,4,12,10\}$. Again we have that $24 \notin S$ [12+12=24], thus it suffices to show that $2$ of the columns $4-12$ contain no elements of $S$. Also $8 \notin S$ [$4+4=8$] and $16 \notin S$ [$4+12=16$]. Once again we have that no element of column $5$ is contained in $S$ and we may thus assume that every other column (among $4-12$) contains an element of $S$. However, we also have that $11 \notin S$ [$1+10=11$] and $22 \notin S$ [$10+12=22$]. Hence no element of column $7$ is contained in $S$, which is a contradiction.

Next we assume that $C=\{1,4,12,20\}$. It follows that $24 \notin S$ [$12+12=24$]; $8 \notin S$ [$4+4=8$]; $16 \notin S$ [$12+4=16$]; $11 \notin S$ [$11+1=12$]; $13 \notin S$ [$12+1=13$]. We record this information in the following table.

\begin{center}
\begin{tabular}{|c|c|c|c|c|c|c|c|c|c|}
  \hline
   4 & 5 & 6 & 7 & 8 & 9 & 10 & 11 & 12 &13 \\ \hline
  7 & \sout{8} & 9 &  \sout{11} &  \sout{13} & 15 & 17 & 25 & 27 & \sout{24}  \\
    14 & \sout{16} & 18 & 22 & 26 & 30 & 34 & 50 & 54 &  \\
   \hline
\end{tabular}
\end{center}

It follows that column $5$ and column $13$ do not contain any element of $S$. Furthermore, since $22+4=26$, at most one of the columns $8$ and $9$ can contain an element of $S$. This completes the proof.

\textbf{Case II: $\{3,12\} \subset S$}. First, $24 \notin S$ [$12+12=24$]. Next, $9 \notin S$ [$9+3=12$] and $15 \notin S$ [$12+3=15$]. Now, $S$ must contain either $18$ or $30$, or columns $6$, $9$ and $13$ would contain no element of $S$ and the proof would be complete. In fact, $S$ must contain exactly one of $18$ and $30$ since $18+12=30$.

We now assume that $30 \in S$ and $18 \notin S$. Then, $54 \in S$ and $27 \notin S$ [$27+3=30$]. If $8 \in S$, then $11 \notin S$ [$3+8=11$] and $22 \notin S$ $[8+22=30]$. This would imply that $S$ contains no elements of columns $6$, $7$, and $13$ and the proof would be complete. We may thus assume that $8 \notin S$ and $16 \in S$. Thus $13 \notin S$ and $26 \in S$ [$13+3=16$]. Also $14 \notin S$ and $7 \in S$ [$14+12=26$]. We now record our progress in the following table

\begin{center}
\begin{tabular}{|c|c|c|c|c|c|c|c|c|}
  \hline
   4 & 5 & 6 & 7 & 8 & 9 & 10 & 11 & 12  \\ \hline
  7 & \sout{8} &\sout{9} &  11 &\sout{13} &\sout{15} & 17 & 25 & \sout{27}  \\
    \sout{14} & 16 & \sout{18} & 22 & 26 & 30 & 34 & 50 & 54   \\
   \hline
\end{tabular}
\end{center}

It follows from this that $C$ cannot contain any of the following elements: $4$ [$4+12=16$], $6$ [$3+3=6$], $5$ [$5+7=12$], and $10$ [$10+16=26$]. If $C$ contains 20, then $17 \notin S$ and $34 \in S$ [$17+3 = 20$]. However, $34+20 = 54$ is then a sum in $S$. Therefore, 20 cannot be in $C$. This would leave at most $3$ elements in $C$ which would contradict our hypothesis.

We now assume that $18 \in S$. We further assume that $16 \in S$. We then have that $30 \notin S$ [$12+18=30$]; $15 \notin S$ [$3+12=15$]; $24 \notin S$ [$12+12=24$]. This implies that columns $9$ and $13$ do not contain elements of $S$. We many now deduce that $17 \in S$ and $34 \notin S$ [$16+18=34$]; $7 \in S$ and $14 \notin S$ [$14+3=17$]; $26 \in S$ and $13 \notin S$ [$13+3=16$]; $22 \in S$ and $11 \notin S$ [$11+7 = 18$]. We again record our progress in a table.

\begin{center}
\begin{tabular}{|c|c|c|c|c|c|c|c|c|c|}
  \hline
   4 & 5 & 6 & 7 & 8 & 9 & 10 & 11 & 12 & 13 \\ \hline
  7 &\sout{8} &\sout{9} &  \sout{11} &\sout{13} &\sout{15} & 17 & 25 & 27 & \sout{24} \\
     \sout{14} & 16 & 18 & 22 & 26 & \sout{30} & \sout{34} & 50 & 54 &   \\
   \hline
\end{tabular}
\end{center}

Thus, $4 \notin C$ [$18+4 = 22$], $1 \notin C$ [$16+1 = 17$], $2 \notin C$ [$16+2 = 18$], and $5 \notin C$ [$5+17 =22$]. This shows that $|C| < 4$, which is a contradiction.

Lastly, we assume that $18 \in S$ and $16 \notin S$. It follows that $24 \notin S $ [$12+12=24$]; $15 \notin S$ [$3+12=15$]; $30 \notin S$ [$12+18=30$]. Hence we may assume that $13 \in S$ and $26 \notin S$ [$18+8=26$]; $50 \in S$ and $25 \notin S$ [$13+12=25$]. We record this below.

\begin{center}
\begin{tabular}{|c|c|c|c|c|c|c|c|c|c|c|}
  \hline
   4 & 5 & 6 & 7 & 8 & 9 & 10 & 11 & 12 & 13  \\ \hline
  7 & 8 & \sout{9} &  11 & 13 & \sout{15} & 17 & \sout{25} & 27 & \sout{24} \\
    14 & \sout{16} & 18 & 22 & \sout{26} & \sout{30} & 34 & 50 & 54 &   \\
   \hline
\end{tabular}
\end{center}

It follows that $C$ does not contain $1$ [$1+12=13$]; $6$ [$6+6=12$]; $5$ [$5+13=18$]; $10$ [$8+10=18$]; $20$ [$8+12=20$]. We now have that $C \subset \{2,4,3,12\}$. But this implies that $|C|\leq 3$ which gives a contradiction.

\textbf{Case III: $\{5,20\} \subset S$}. Using the fact that $\{1,4\} \nsubseteq C$ and $\{3,12\} \nsubseteq C$, we have that $C$ either contains $\{4,5,20\}$ or $\{12,5,20\}$, or is equal to $\{1,3,5,20\}$ or $\{2,6,5,20\}$.

If $C$ contains either $\{4,5,20\}$ or $\{12,5,20\}$, we must have that $24 \notin S$ [$4+20=24$, $12+12=24$], thus it suffices to show that $2$ of the columns $4-12$ contain no elements of $S$.

Assume that $\{4,5,20\} \subset C$. Hence $8 \notin S$ [$4+4=8$] and $16 \notin S$ [$4+16=20$]. In addition, $15 \notin S$ [$15+5=20$] and $25 \notin S$ [$5+20=25$]. Recording this in a table we have that

\begin{center}
\begin{tabular}{|c|c|c|c|c|c|c|c|c|c|}
  \hline
   4 & 5 & 6 & 7 & 8 & 9 & 10 & 11 & 12  \\ \hline
  7 &\sout{8} & 9 &  11 &13 &\sout{15} & 17 &\sout{25} & 27  \\
      14 &\sout{16} & 18 & 22 & 26 & 30 & 34 & 50 & 54   \\
   \hline
\end{tabular}
\end{center}

Now $S$ cannot contain an element of column $5$. Moreover $S$ can only contain an element from either column $9$ or $11$ since $20+30=50$. Hence there are two columns (among $4-12$) that do not contain an element of $S$ and the proof is complete.

Next, assume that $\{12,5,20\} \subset C$. It follows that $15 \notin S$ [$5+15=20$], $25 \notin S$ [$20+5 = 25$], $7 \notin S$ [$5+7=12$], $8 \notin S$ [$8+12=20$], $17 \notin S$ [$12+5=17$]. We record this information in a table.

\begin{center}
\begin{tabular}{|c|c|c|c|c|c|c|c|c|}
  \hline
   4 & 5 & 6 & 7 & 8 & 9 & 10 & 11 & 12  \\ \hline
  \sout{7} & \sout{8} & 9 &  11 &13 &\sout{15} &\sout{17} &\sout{25} & 27  \\
      14 & 16 & 18 & 22 & 26 & 30 & 34 & 50 & 54   \\
   \hline
\end{tabular}
\end{center}

Now notice that only one element from $\{30,50\}$ can be in $S$ [$20+30=50$]. Furthermore, only element from $\{14,34\}$ can be in $S$ [$14+20=34$]. Thus two columns (among $4-12$) must not contain elements of $S$ and the subcase is complete.

Next we assume that $C=\{1,3,5,20\}$. We first suppose that $24 \in S$. It follows that $8 \notin S$ [$3+5=8$]; $15 \notin S$ [$15+5=20$]; $17 \notin S$ [$17+3=20$]; $25 \notin S$ [$20+5=25$]; $27 \notin S$ [$3+24=7$]. We record this below.

\begin{center}
\begin{tabular}{|c|c|c|c|c|c|c|c|c|c|}
  \hline
   4 & 5 & 6 & 7 & 8 & 9 & 10 & 11 & 12 & 13  \\ \hline
  7 & \sout{8} & 9 &  11 & 13 & \sout{15} & \sout{17} & \sout{25} & \sout{27} & 24  \\
      14 & 16 & 18 & 22 & 26 & 30 & 34 & 50 & 54 &   \\
   \hline
\end{tabular}
\end{center}

If follows that only one of the columns $9$ and $11$ can contain an element of $S$ [$30+20=50$]. Furthermore, only one of the columns $10$ and $12$ can contain an element of $S$. We may then assume that $16 \in S$. Thus $26 \in S$ and $13 \notin S$ [$13+3=16$]; $22 \in S$ and $11 \notin S$ [$11+5=16$]. From this we conclude that $50 \notin S$ [$26+24=50$], hence $30 \in S$. But now $54 \notin S$ [$24+30=54$] and $34$ must be in $S$. This implies that $7 \in S$ and $14 \notin S$ [$14+16=30$]; $9 \in S$ and $18 \notin S$ [$18+16=34$]. However, we have reached a contradiction since $9+7=16$.

We now assume that $C=\{1,3,5,20\}$ and $24 \notin S$. Again we have that $15 \notin S$ [$15+5=20$]; $25 \notin S$ [$5+20=25$]; $17 \notin S$ [$17+3=20$]. We record this below.

\begin{center}
\begin{tabular}{|c|c|c|c|c|c|c|c|c|c|}
  \hline
   4 & 5 & 6 & 7 & 8 & 9 & 10 & 11 & 12 & 13  \\ \hline
  7 & 8 & 9 &  11 & 13 & \sout{15} & \sout{17} & \sout{25} & 27 & \sout{24}  \\
      14 & 16 & 18 & 22 & 26 & 30 & 34 & 50 & 54 &   \\
   \hline
\end{tabular}
\end{center}

Again, it follows that only one of the columns $9$ and $11$ can contain an element of $S$ [$20+30=50$]. If $27 \notin S$, then only one of the columns $10$ and $12$ can contain an element of $S$, and the proof would be complete. We may now assume that $27 \in S$. This implies that $7 \notin S$ and $14 \in S$ [$7+20=27$]. However, we now have that $34 \notin S$ [$20+14=34$] which implies that column $10$ does not contain an element of $S$, and the proof of this case is  complete.

Lastly, we consider $C=\{2,5,6,20\}$. We then have that $15 \notin S$ [$15+5=20$]; $25 \notin S$ [$20+5=25$]; $11 \notin S$ [$5+6=11$]; $22 \notin S$ [$20+2=22$]; $18 \notin S$ [$18+2=20$]; $8 \notin S$ [$2+6=8$]; $14 \notin S$ [$6+14=20$]. We record this in the table below.

\begin{center}
\begin{tabular}{|c|c|c|c|c|c|c|c|c|}
  \hline
   4 & 5 & 6 & 7 & 8 & 9 & 10 & 11 & 12  \\ \hline
  7 & \sout{8} & 9 &  \sout{11} &13 & \sout{15} & 17 & \sout{25} & 27  \\
      \sout{14} & 16 & \sout{18} & \sout{22} & 26 & 30 & 34 & 50 & 54   \\
   \hline
\end{tabular}
\end{center}

It follows that no element of column $7$ is contained in $S$. Furthermore, either column $9$ or $11$ must not contain an element of $S$. Lastly one of the columns $4-6$ must not contain an element of $S$ since $7+9=16$. This completes the proof.

\section{Case 4: $|C| = 3$ and $24 \notin S$}

To show that $S$ cannot be a sum-free set of size $12$, it suffices to show that two of the columns $4-13$ do not contain elements of $S$. Here we assume that $24 \notin S$. Thus it suffices to show that one of the columns $4-12$ does not contain an element of $S$. We suppose this is false, which implies that exactly one element of $\{7,14\}$ and one element of $\{8,16\}$ are contained in $S$. There are $4$ possible combinations of such elements and we consider each of these four cases separately.

\textbf{Case I: $7,8 \in S$}. We have that $30 \in S$ and $15 \notin S$ [$7+8=15$]. Also $11 \in S$ and $22 \notin S$ [$8+22=30$]. Next $9 \in S$ and $18 \notin S$ [$7+11=18$]. We can conclude that $1 \notin C$ [$1+7=8$], $2\notin C$ [$2+7=9$], $4 \notin C$ [$4+7=11$], $3 \notin C$ [$3+8=11$], and $20 \notin C$ [$9+11=20$]. We record this in the following table.

\begin{center}
\begin{tabular}{|c|c|c|}
  \hline
  1 &  2 &  3  \\ \hline
  \sout{1} & \sout{3} &  5  \\
  \sout{2} & 6 & 10   \\
  \sout{4} & 12 & \sout{20}   \\
  \hline
\end{tabular}
\end{center}

Now $C$ can have size at most two which gives a contradiction.

\textbf{Case II: $7,16 \in S$}. We have that $18 \in S$ and $9 \notin S$ [$9+7=16$]. Furthermore, $50 \in S$ and $25 \notin S$ [$7+18=25$]; $17 \in S$ and $34 \notin S$ [$16+18=34$]; $22 \in S$ and $11 \notin S$ [$11+7=18$]; $30 \in S$ and $15 \notin S$ [$15+7=22$]; $26 \in S$ and $13 \notin S$ [$13+17=30$].

We now have that $1\notin C$ [$1+16=17$], $2 \notin C$ [$2+16=18$], $4 \notin C$ [$4+18=22$], $6 \notin C$ [$6+16=22$], $12 \notin C$ [$12+18=30$], $10 \notin C$ [$10+7=17$], and $20 \notin C$ [$20+30=50$]. This implies that $|C|\leq 2$ and gives a contradiction.

\textbf{Case III: $14,8 \in S$}. We have $11 \in S$ and $22 \notin S$ [$8+14=22$]. Now $3 \notin C$ [$3+8=11$], $4 \notin C$ [$4+4=8$] and $6 \notin C$ [$6+8=14$]. Now suppose that $12 \in C$. We would then have that $1 \notin C$ [$1+11=12$], $2 \notin C$ [$2+12=14$] and $20 \notin S$ [$8+12=20$]. As the following table indicates, this implies that $|C|\leq 2$ and gives a contradiction.

\begin{center}
\begin{tabular}{|c|c|c|}
  \hline
  1 &  2 &  3  \\ \hline
  \sout{1} & \sout{3} & 5  \\
  \sout{2} & \sout{6} & 10   \\
  \sout{4} & 12 & \sout{20}   \\
  \hline
\end{tabular}
\end{center}

We may now assume that $\{3,6,12\} \cap C = \emptyset$. Since either column $1$ or $3$ must now contain two elements of $C$, it follows that either $\{1,4\} \subset C$ or $\{5,20\} \subset C$. However, we have previously shown that $4 \notin C$ [$4+4=8$], hence we have that $\{5,20\} \subset C$. This implies that $15 \notin S$ [$15+5=20$], $25 \notin S$ [$20+5 = 25$] and $50 \in S$, so $30 \notin S$ [$20+30=50$]. Hence, $S$ contains no element of column $9$, which gives a contradiction.

\textbf{Case IV: $14,16 \in S$}. First we assume that $12 \in S$ as well. It follows that $15 \in S$ and $30 \notin S$ [$14+16=30$]. Also $13 \in S$ and $26 \notin S$ [$14+12=26$]. Next, $50 \in S$ and $25 \notin S$ [$12+13=25$]. Furthermore, $54 \in S$ and $27 \notin S$ [$12+15=27$]. We record this information in the following table.

\begin{center}
\begin{tabular}{|c|c|c|c|c|c|c|c|c|}
  \hline
   4 & 5 & 6 & 7 & 8 & 9 & 10 & 11 & 12  \\ \hline
  \sout{7} & \sout{8} & 9 &  11 & 13 & 15 & 17 & \sout{25} & \sout{27}  \\
      14 & 16 & 18 & 22 & \sout{26} & \sout{30} & 34 & 50 & 54   \\
   \hline
\end{tabular}
\end{center}

Now this implies that $1 \notin C$ [$1+13=14$], $2 \notin C$ [$2+14=16$], $4 \notin C$ [$4+50=54$], $3 \notin C$ [$3+13=16$] and $6 \notin C$ [$6+6=12$]. We record this in the following table

\begin{center}
\begin{tabular}{|c|c|c|}
  \hline
  1 &  2 &  3  \\ \hline
  \sout{1} & \sout{3} &  5  \\
  \sout{2} & \sout{6} & 10   \\
  \sout{4} & 12 & 20   \\
  \hline
\end{tabular}
\end{center}

It follows that $5$ and $20$ must both be in $C$, however this gives us a contradiction since $5+15=20$.

Now we assume that $12 \notin S$. As before, we may assume that $15 \in S$ and $30 \notin S$ [$14+16=30$]. However, this implies that $1 \notin C$ [$1+14=15$] and $2\notin C$ [$2+14=16$], thus we have:

\begin{center}
\begin{tabular}{|c|c|c|}
  \hline
  1 &  2 &  3  \\ \hline
  \sout{1} &  3 &  5  \\
  \sout{2} & 6 & 10   \\
   4 & \sout{12} & 20   \\
  \hline
\end{tabular}
\end{center}

We claim that $20 \notin C$ as well. If $20 \in C$ then $5 \notin C$ [$5+15=20$], and thus $4 \in C$ (since $C$ would then contain at most one element from each of columns $2$ and $3$). However, this would give a contradiction since $4+16=20$. Thus $20 \notin C$. Now $C$ can contain at most one element from each of the columns $1-3$, and thus it must contain $4$. Hence $22 \in S$ and $11 \notin S$ [$4+11=15$]. Furthermore, $9 \in S$ and $18 \notin S$ [$4+18=22$]. Also $13 \in S$ and $26 \notin S$ [$4+22=26$]. Thus $6 \notin C$ [$6+16=22$] and $5 \notin C$ [$5+9=14$]. It follows that $C=\{4,3,10\}$, however this contradicts our claim that $13 \in S$.

\section{Case 5: $|C| = 3$ and $24 \in S$}We assume that $ 24 \in S$. Furthermore, let us assume that none of the elements $\{2,3,6,10\}$ are contained in $S$. We record this information with the fact that $12 \notin C$ [$12+12=24$] in the following table.

\begin{center}
\begin{tabular}{|c|c|c|}
  \hline
  1 &  2 &  3  \\ \hline
  1 &  \sout{3} &  5  \\
  \sout{2} &  \sout{6} & \sout{10}   \\
   4 & \sout{12} & 20   \\
  \hline
\end{tabular}
\end{center}

It follows that either $\{1,4\} \subset C$ or $\{5,20\}\subset C$. In the first case, $\{1,4\} \subset C$, we have that $5 \notin C$ [$1+4=5$] and $20 \notin C$ [$4+20=24$]. Thus we may assume that $\{5,20\}\subset C$. Moreover, $4 \notin C$ [$4+20=24$], thus $C=\{1,5,20\}$. Assuming that $C=\{1,5,20\}$ and $24 \in S$, we have that $15 \notin S$ [$15+5=20$] and $25 \notin S$ [$20+5=25$]. We record this in the table below.

\begin{center}
\begin{tabular}{|c|c|c|c|c|c|c|c|c|c|}
  \hline
   4 & 5 & 6 & 7 & 8 & 9 & 10 & 11 & 12 & 13 \\ \hline
   7 & 8 & 9 &  11 & 13 & \sout{15} & 17 & \sout{25} & 27 & 24  \\
      14 & 16 & 18 & 22 & 26 & 30 & 34 & 50 & 54 &  \\
   \hline
\end{tabular}
\end{center}

Since $20 \in C$, only one of the integers $30$ and $50$ can be contained in $S$. Thus it suffices to show that if $30 \in S$ (respectively $50 \in S$) there is a column (among $4-12$) in addition to $11$ (respectively $9$) that contains no element of $S$.

If $30 \in S$, we then have that $27 \in S$ and $54 \notin S$ [$30+24=54$]. Next $14 \in S$ and $7 \notin S$ [$20+7=27$]. Furthermore, $26 \in S$ and $13 \notin S$ [$13+1 =14$]. This gives a contradiction since $26+1=27$.

Next, if $50 \in S$ we have that $13 \in S$ and $26 \notin S$ [$24+26=50$]. Also $22 \in S$ and $11 \notin S$ [$11+13=24$]. Furthermore, $54 \in S$ and $27 \notin S$ [$22+5=27$]. Next $17 \in S$ and $34 \notin S$ [$34+20=54$]. However, we have reached a contradiction since $17+5=22$.

It now suffices to show that we may assume that $2,3,6,10 \notin S$.

\subsection{$2 \notin S$}

If $2 \in S$, we have that $22 \notin S$ [$2+22=24$] and $26 \notin S$ [$2+24=26$]. This implies that exactly one of $11$ and $13$ are in $S$. We consider these two cases separately.

If $13 \in S$, we have that $30 \in S$ and $15 \notin S$ [$13+2=15$]. Next $34 \in S$ and $17 \notin S$ [$13+17=30$]. Also $14 \in S$ and $7 \notin S$ [$7+24=34$]. Lastly $27 \in S$ and $54 \notin S$ [$24+30=54$]. From this we can conclude that $1 \notin C$ [$1+13=14$], $4 \notin C$ [$4+30=34$], $3 \notin C$ [$3+24=27$], $6 \notin C$ [$6+24=30$], $12 \notin C$ [$12+12=24$], $10 \notin C$ [$10+24=34$], and $20 \notin C$ [$20+14=34$]. It follows that $|C| \leq 2$, which gives a contradiction.

If $11 \in S$, we have $22, 13, 26 \notin S$, $18 \in S$ and $9\notin S$ [$2+9=11$]. Next $14 \in S$ and $7 \notin S$ [$7+11=18$]. Also $8 \in S$ and $16 \notin S$ [$14+2=16$]. Furthermore, $50 \in S$ and $25 \notin S$ [$11+14=25$]. We now have that $1 \notin C$ [$1+1=2$], $4 \notin C$ [$2+2=4$], $3 \notin C$ [$8+3=11$], $6 \notin C$ [$18+6=24$], $12 \notin C$ [$12+12=24$], $10 \notin C$ [$14+10=24$], and $20 \notin C$ [$2+18=20$]. It follows that $|C| \leq 2$, which gives a contradiction.

\subsection{$10 \notin S$}

Assuming $10,24 \in S$, we have that neither $14$ or $34$ are contained in $S$. Hence exactly one element of the set $\{7,17\}$ must be contained in $S$ [$7+17=24$]. Furthermore, we may assume that $6 \notin C$, since $6 \in C$ would imply that $18 \notin S$ [$18+6=24$], $30 \notin S$ [$24+6=30$]. In addition, since $9+15=24$, $S$ can only contain one element of the set $\{9,15\}$.  This implies that either column $6$ or column $9$ does not contain an element of $S$ which contradicts our hypothesis.

We consider the two cases mentioned above separately. First, assuming that $7 \in S$, we record the implications in the following table.

\begin{center}
\begin{tabular}{|c|c|c|c|c|c|c|c|c|c|}
  \hline
   4 & 5 & 6 & 7 & 8 & 9 & 10 & 11 & 12 & 13 \\ \hline
   7 & 8 & 9 &  11 & 13 & 15 & \sout{17} & 25 & 27 & 24  \\
     \sout{14} & 16 & 18 & 22 & 26 & 30 & \sout{34} & 50 & 54 &  \\
   \hline
\end{tabular}
\end{center}

We have $2,6,12 \notin C$ by previous considerations ([$12+12 = 24$], $6 \notin C$ from above, and $2\notin C$ from the previous subsection). Furthermore, we have that $5\notin C$ [$5+5=10$] and $20 \notin C$ [$10+10=20$]. Lastly we have that $3 \notin C$ [$3+7=10$]. This implies that $C=\{1,4,10\}$. However, we now have that $18 \in S$ and $9 \notin S$ [$9+1=10$]. In addition $22 \in S$ and $11 \notin S$ [$10+1=11$]. However, this gives a contradiction since $18+4=22$.

Secondly we assume that $17 \in S$. We record this in tabular form.

\begin{center}
\begin{tabular}{|c|c|c|c|c|c|c|c|c|c|}
  \hline
   4 & 5 & 6 & 7 & 8 & 9 & 10 & 11 & 12 & 13 \\ \hline
   \sout{7} & 8 & 9 &  11 & 13 & 15 & 17 & 25 & 27 & 24  \\
    \sout{14} & 16 & 18 & 22 & 26 & 30 & \sout{34} & 50 & 54 &  \\
   \hline
\end{tabular}
\end{center}

We deduce that $54 \in S$ and $27 \notin S$ [$17+10=27$]. Also $15 \in S$ and $30 \notin S$ [$30+24=54$]. Furthermore, $50 \in S$ and $25 \notin S$ [$15+10=25$]. Lastly $13 \in S$ and $26 \notin S$ [$24+26=50$]. We now have that $2,6,12,5,20 \notin C$ by previous considerations. This implies that $C \subset \{1,3,4,10\}$. However, $4 \notin C$ [$4+13=17$] and $3 \notin C$ [$3+10=13$], which gives a contradiction.

\subsection{$6 \notin S$}

We may assume that $12 \notin C$ [$12+12=24$] and that $2,10 \notin C$ (from the previous cases). If $6 \in S$, we have that $18 \notin S$ [$18+6=24$] and $30 \notin S$ [$24+6=30$]. This implies that $S$ must contain exactly one element of $\{9,15\}$, and one of the columns $6$ or $9$ cannot contain an element of $S$.

Assume that $15 \in S$. It follows that $1 \notin C$, since $1 \in C$ would imply that $14,16 \notin S$ [$1+14=15$, $1+15=16$]. Furthermore, we would have that at most one of the elements $7,8$ is contained in $S$ and, hence, at most one of the columns $4,5$ contains an element of $S$. This implies that $C \subset \{4,5,6,20\}$. However, it follows that $20 \notin C$, since $20 \in C$ would imply that $C=\{4,6,20\}$ or $\{5,6,20\}$ (recall we are assuming $6 \in C$), but $15+5=20$ and $4+20=24$. Thus $C=\{4,5,6\}$. However, this implies that $22 \in S$ and $11 \notin S$ [$11+4=15$]; $13 \in S$ and $26 \notin S$ [$22+4=26$]; $8 \in S$ and $16 \notin S$ [$16+6=22$]. However, this gives a contradiction since $4+4=8$.

We now assume that $9 \in S$. It follows that $20 \notin C$, since $20 \in C$ would imply that $7 \in S$ and $14 \notin S$ [$14+6=20$]. Also, $34 \in S$ and $17 \notin S$ [$17+7=24$];  $54 \in S$ and $27 \notin S$ [$27+7=34$]. However, this gives a contradiction since $34+20=54$.

It follows that $C \subset \{1,4,5,6\}$, and $6 \in C$. Also $C \neq \{1,5,6\}$ [$1+5=6$] and $C \neq \{4,5,6\}$ [$4+5=9$]. It follows that $C=\{1,4,6\}$. We then have that $26 \in S$ and $13 \notin S$ [$9+4=13$]; $50 \in S$ and $25 \notin S$ [$1+24=25$]. We have reached a contradiction however since $26+24=50$.

\subsection{$3 \notin S$} We collect the results from the preceding three subsections in the following table.

\begin{center}
\begin{tabular}{|c|c|c|}
  \hline
  1 &  2 &  3  \\ \hline
  1 &  3 &  5  \\
  \sout{2} &  \sout{6} & \sout{10}   \\
   4 & \sout{12} & 20   \\
  \hline
\end{tabular}
\end{center}

If $3 \in C$, it follows that $C$ must be one of the following four sets: $\{3,5,20\}$, $\{3,4,5\}$, $\{1,3,5\}$, $\{1,3,20\}$. We note that $\{4,3,20\}$ is not possible because $24 \in S$ [$4+20 = 24$].

If $C=\{3,5,20\}$, then $27 \notin S$ [$3+24=27$], $25 \notin S$ [$5+20=25$], $17 \notin S$ [$3+17=20$], $15 \notin S$ [$15+5=20$]. We record these results in the following table.

\begin{center}
\begin{tabular}{|c|c|c|c|c|c|c|c|c|c|}
  \hline
   4 & 5 & 6 & 7 & 8 & 9 & 10 & 11 & 12 & 13 \\ \hline
   7 & 8 & 9 &  11 & 13 & \sout{15} & \sout{17} & \sout{25} & \sout{27} & 24  \\
      14 & 16 & 18 & 22 & 26 & 30 & 34 & 50 & 54 &  \\
   \hline
\end{tabular}
\end{center}

Notice that $S$ can include only one element of the following two disjoint sets $\{30,50\}$ and $\{34,54\}$. It follows that two columns in the table above will not contain an element of $S$, which yields a contradiction.

If $C=\{3,4,5\}$, we have that $7 \notin S$ [$3+4=7$]; $9 \notin S$ [$5+4=9$]; $8 \notin S$ [$4+4=8$]. We record this information in the following table.

\begin{center}
\begin{tabular}{|c|c|c|c|c|c|c|c|c|c|}
  \hline
   4 & 5 & 6 & 7 & 8 & 9 & 10 & 11 & 12 & 13 \\ \hline
   \sout{7} & \sout{8} & \sout{9} &  11 & 13 & 15 & 17 & 25 & 27 & 24  \\
      14 & 16 & 18 & 22 & 26 & 30 & 34 & 50 & 54 &  \\
   \hline
\end{tabular}
\end{center}

We have that either column $4$ or $6$ must not contain an element of $S$ [$14+4=18$]. Since the proof is complete once we show that two columns among $4-12$ do not contain an element of $S$, it follows that $16 \in S$; $26 \in S$ and $13 \notin S$ [$13+3=16$]; $11 \in S$ and $22 \notin S$ [$22+4=26$]. This give a contradiction, however, since $11+5=16$.

If $C=\{1,3,20\}$, we have that $17 \notin S$ [17+3=20]; $25 \notin S$ [$1+24=25$]; $27 \notin S$ [$24+3=27$]. Notice that $S$ must contain exactly one element from columns $10$ and $12$, either $34$ or $54$ [34+20=54]. If $34 \in S$, we have that $7 \in S$ and $14 \notin S$ [$14+20=34$]; $26 \in S$ and $13 \notin S$ [$13+7=20$]. This implies that $50 \notin S$ [$24+26 =50$] and that column $11$ does not contain an element of $S$, a contradiction.

If $C=\{1,3,5\}$, we have that $27 \notin S$ [$24+3=27$]; $25 \notin S$ [$1+24=25$]; $8 \notin S$ [$3+5=8$]. We first assume that $16 \in S$. It follows that either column $10$ or $11$ does not contain an element of $S$ since $17 \notin S$ [$16+1=17$] and $34+16=50$. Thus $54 \in S$, $15 \in S$ and $30 \notin S$ [$24+16=30$]. However this gives a contradiction since $15+1=16$. Thus it suffices to assume that $16 \notin S$, and $S$ does not contain an element of column $5$. It follows that $50 \in S$ and $54 \in S$. In addition $15 \in S$ and $30 \notin S$ [$24+30=54$]; $9 \in S$ and $18 \notin S$ [$15+3=18$]. However, we have reached a contradiction since $15+9=24$.

This completes the proof.

\bigskip
\hrule
\bigskip

\noindent 2010 {\it Mathematics Subject Classification}:
Primary 11B30.

\noindent \emph{Keywords: } sum-free.

\bigskip

\vskip .1in

\end{document}